\documentclass[12pt]{article}
\usepackage{amssymb,latexsym,amsmath,graphicx} 
\usepackage{booktabs}
\usepackage{enumitem}
\usepackage{latexsym}
\usepackage{mathrsfs}
\usepackage{float}
\usepackage{lscape}
\usepackage[mathscr]{euscript}
\newtheorem{thm}{Theorem}

\newtheorem{rem}{Remark}

\begin{document}
\begin{center}
\Large{{\bf A simple approach to construct confidence bands for a regression function with incomplete data}}
\end{center}

 \vskip .1 in

 \begin{center}
{\bf Ali Al-Sharadqah}\footnote{Email: ali.alsharadqah@csun.edu}
{\bf and Majid Mojirsheibani}\footnote{Corresponding author. ~Email:  majid.mojirsheibani@csun.edu ~(Tel.: 1-818-677-7814)
	
This work is supported by the NSF Grant DMS-1407400
of Majid Mojirsheibani.}

{\bf Department of Mathematics, California State University Northridge, CA, 91330, USA}\\

\end{center}

\begin{abstract}
A long-standing problem in the construction of asymptotically correct confidence bands for a regression function $m(x)=E[Y|X=x]$, where $Y$ is the response variable influenced by the covariate $X$, involves the situation where $Y$ values may be missing at random, and where the selection probability, the density function $f(x)$ of $X$, and the conditional variance of $Y$ given $X$ are all completely unknown. This can be particularly more complicated in nonparametric situations. In this paper we propose a new kernel-type regression estimator and study the limiting distribution of the properly normalized versions of the maximal deviation of the proposed estimator from the true regression curve. The resulting limiting distribution will be used to construct uniform confidence bands for the underlying regression curve with asymptotically correct coverages. The focus of the current paper is on the case where $X\in \mathbb{R}$. We also perform numerical studies to assess the finite-sample performance of the proposed method. In this paper, both mechanics and the theoretical validity of our methods are discussed.
\end{abstract}


\allowdisplaybreaks
\vspace{8mm}
\noindent
\section{Introduction}
Nonparametric regression estimation has important applications in both statistical estimation theory and classical statistical inferential procedures such as tests of hypotheses or the construction of uniform confidence bands for a true regression function.  Confidence bands, in particular,  provide insight into the variability of the estimators of the entire regression curve and can also be used to study and investigate certain global features, such as the shape, of the true curve. A very standard procedure to construct asymptotically correct confidence bands for a regression function, over a connected compact set, is based on the limiting distribution of the properly normalized versions of the maximal deviation of the regression estimator from the true regression curve. 
Results along these lines include the work of Johnston (1982) who constructed confidence bands based on kernel regression estimators; H\"ardle (1989) established uniform
confidence bands for M-smoothers; Eubank and Speckman (1993) proposed bias-corrected confidence bands based on nonparametric kernal regression with fixed design points, whereas  Xia (1998) considered random design points under dependence. Additionally, bootstrap confidence bands  based on nonparametric regression have been proposed by  Neumann and Polzehl (1998), Claeskens and Keilegom (2003), and Song et al. (2012). H\"ardle and Song (2010) constructed uniform confidence bands for a quantile regression curve with a one-dimensional predictor, whereas  Cai et al. (2014) constructed adaptive confidence bands based on nonparametric regression functions. Mass\'{e} and Meiniel (2014) developed adaptive confidence bands for the case of nonparametric fixed design regression models, and Proksch (2016) developed uniform confidence bands in a nonparametric regression setting with deterministic and multivariate predictor. Another related result is that of Gu and Yang (2015).

\vspace{4mm}\noindent
The papers cited above as well as most of the results in the literature deal with the cases where the data are fully observable.  
The focus of this paper is on the realistic case where the response variable could be unobservable or {\it missing}. More specifically, let $(X,Y)$ be a random pair with the cumulative distribution function (cdf) $G(x,y)$.  Here, $X$ is a $d$-dimensional random vector of covariates and $Y$ is the response variable influenced by $X$. Given the independent and identically distributed (iid) data  $(X_1, Y_1), \dots , (X_n, Y_n)$ from $G(x,y)$, let $m_n(x)$ be the Nadaraya-Watson (Nadaraya (1970), Watson (1964)) kernel regression estimator of the regression function $m(x):=E[Y| X= x]$, i.e., 
\begin{equation}\label{D1}
m_n(x) =  \sum_{i=1}^n Y_i K((x-X_i)/h_n) \div \sum_{i=1}^n K((x-X_i)/h_n)\,,
\end{equation}
where 0/0 :=\,0 by convention. Here $K:\mathbb{R}^d\to\mathbb{R}_+$ is the kernel used with the {\it bandwidth} $h_n$. However, our focus in this paper is on the case of $d=1$. Now, for various reasons, some of the $Y_i$'s may be unavailable or missing from the data. Missing data are common in opinion polls, survey data, mail questionnaires, data collected in medical research and other scientific studies. In this paper we consider the case where $Y$ may be Missing At Random (MAR). More specifically, let $\Delta$ be the Bernoulli random variable defined as $\Delta=0$ if $Y$ is missing and $\Delta=1$, otherwise. Then, the MAR assumption states:
\begin{equation} \label{MAR}
P\{\Delta=1|X=x, Y=y\}= P\{\Delta=1|X=x\} = E[\Delta|X=x] ~=:~ p(x),
\end{equation}
i.e., the probability that $Y$ is missing does not depend on $Y$ itself. For more on this and other types of missing probability mechanism see, for example, Little and Rubin (2002). Here, the missing probability mechanism $p(x)$, also called the {\it selection probability}, is assumed to be completely unknown.  In the rest of this paper, the iid data will be represented by $(X_1, Y_1, \Delta_1),\dots,(X_n, Y_n, \Delta_n)$. Some work has been done for the simpler problem of constructing confidence intervals for $m(x_0)=E[Y|X=x_0]$, where $x_0$ is a given point in $\mathbb{R}^d$; see, for example, Qin et al. (2014) as well as Lei and Qin (2011).  However, to the best of our knowledge, a long-standing problem in constructing uniform confidence bands for the regression function $m(x)$, over compact sets, involves the situation where the response variable $Y$ may be missing at random and the function $p(x)$ in (\ref{MAR}), the density function $f(x)$ of $X$, and the conditional variance $\sigma_0^2(x)=E[Y^2|X=x]-m^2(x)$ are all completely unknown. Of course, it should be possible to form asymptotically correct uniform confidence bands for $m(x)$ under the restrictive assumptions that $p(x)$, $f(x)$, and $\sigma_0^2(x)$, (or certain functions of these quantities) are known.  However, since such assumptions are unrealistic and not warranted in practice, they will not be pursued in this paper. \\

\noindent
In passing we also note that in the case of censored data, Hollander et al. (1997)  proposed confidence bands for survival functions based on the  empirical likelihood method. Li and van Keilegom (2002) constructed confidence bands for the conditional survival function under random censorship. Wang and Shen (2008) constructed confidence bands of a conditional survival function when the censoring indicators are missing at random. Mondal and Subramanian (2016) developed simultaneous confidence bands for Cox regression in a semiparametric random censorship setup. In another closely related result, Wang and Qin (2010) constructed empirical likelihood confidence bands for a distribution function with missing responses. \\

\noindent
In the next section we propose a new kernel-type regression estimator with missing response variables. Our main result in Theorem \ref{THM-2} deals with the limiting distribution of the properly normalized versions of the maximal deviation of the proposed estimator from the true regression curve. Theorem \ref{THM-2} may be viewed as a counterpart of the classical result of Liero (1982) for the setup with no missing data; see Theorem \ref{THM-1} in Section \ref{sub1}. Our results will be used to develop a new effective procedure for constructing uniform confidence bands for a regression function in the presence of missing response variables with asymptotically correct coverages. Our numerical results also confirm the finite-sample effectiveness of our procedures.

\section{Main results}
\subsection{Preliminaries and the background tools}\label{sub1}
To provide the necessary background tools, let $m_n(x)$ be the kernel regression estimator defined in (\ref{D1}). Also, let 
\begin{equation} \label{sig0}
\widehat{f}_n(x)=\frac{1}{nh_n}\sum_{i=1}^n K((x-X_i)/h_n)\,~~\mbox{and}~~~
\sigma^2_n(x) = \frac{\sum_{i=1}^n Y^2_i K((x-X_i)/h_n)}{\sum_{i=1}^n K((x-X_i)/h_n)} - m^2_n(x)
\end{equation}
be, respectively, the kernel estimators of the density $f$ of $X$ and the conditional variance $\sigma_0^2(x)=E[Y^2|X=x]-m^2(x)$. The limiting distribution of the properly normalized versions of the statistic \, $\sup_{x\in[0,1]}\,\sqrt{\widehat{f}_n(x)/\sigma^2_n(x)}~\big| m_n(x)-m(x)\big|$ have been studied by many authors; see, for example Wandl (1980), Johnston (1982), Liero (1982), and H\"ardle (1990). In passing, we also note that the interval $[0,1]$ may be any connected compact subset of the interior of the support of $f$.
For the case where $X$ is a $d$-dimensional vector, one may refer to the results of Konakov and Piterbarg (1984) and those of Muminov (2011, 2012). 
To state our proposed estimators and results, we first state a number of classical assumptions, some of which will also be used throughout this paper. These assumptions are virtually all the same as those in Liero (1982). 

\vspace{2mm}\noindent
{\bf Assumption (A).} The random pair $(X, Y)$ has a  probability density function (pdf), $g(x,y)$, with respect to the Lebesgue measure. The random variable $Y$ is almost surely bounded, i.e., $P\{B_1 \leq Y\ \leq B_2\}=1$ for constants $-\infty < B_1 <B_2 <\infty$.

\vspace{3mm}\noindent
{\bf Assumption (B).} The pdf of $X$, $f$, is strictly positive on $[0,1]$ and vanishes outside of a finite interval $[a,b]$, where $[0,1]\subset (a,b)$.

\vspace{3mm}\noindent
{\bf Assumption (C).} The functions $m(x)$, $f(x)$, and $\sigma_0^2(x)=E[(Y-m(X))^2|X=x]$ are twice differentiable with bounded derivatives. Furthermore, $\sigma_0^2(x)$ is strictly positive on $[0,1]$.

\vspace{3mm}\noindent
To state the next assumption, put $Z=Y-m(X)$ and let $\widetilde{G}(x,z)$ and $\widetilde{g}(x,z)$ be the cdf and the pdf of the vector $(X,Z)$. Also, let $F$ be the cdf of $X$ and define $H(z|x)$ and $h(z|x)$ to be the conditional cdf and the conditional pdf of $Z$ given $X$, respectively. 

\vspace{3mm}\noindent
{\bf Assumption (D).} ~$\widetilde{g}^{1/2}(x,z)$ is differentiable with respect to both $x$ and $z$, and the partial derivatives are bounded. Furthermore, the inverse functions $H^{-1}$ and $F^{-1}$ of $H$ and $F$ exist and $\frac{\partial}{\partial x} H^{-1}\big(z|F^{-1}(x)\big)$ and $\frac{\partial}{\partial z} H^{-1}\big(z|F^{-1}(x)\big)$ are bounded.

\vspace{3mm}\noindent
{\bf Assumption (E).} The kernel $K$ is a density function and has a bounded support $[-A, A]$ for some $A>0$. Furthermore, $K$ is continuously differentiable and satisfies $\int x \, K(x)\, dx =0$.

\vspace{5mm}\noindent
We have the following classical result (see, for example, Liero (1982)). 

\begin{thm} \label{THM-1}
	Let $h_n=n^{-\delta}$, $\frac{1}{5} < \delta < \frac{1}{3}$, and suppose that assumptions (A)-(E) hold. Then
	\begin{equation}\label{D2}
	P\left\{
	\sqrt{2\delta\log n}\left(
	\sqrt{\frac{nh_n}{c_K}} \,\, \sup_{x\in[0,1]}\,\sqrt{\frac{\widehat{f}_n(x)}{\sigma^2_n(x)}}~\Big| m_n(x)-m(x)\Big|-d_n\right)\,\leq\, y\right\}\rightarrow  \exp\left(-2\,e^{-y}\right),
	\end{equation}
	as $n\to\infty$, where $c_K=\int K^2(t)\,dt$ and  
	\begin{eqnarray}
	d_n &=& \sqrt{2\delta\log n} \,\,+\,
	\left\{\begin{array}{ll}
	\frac{\log (C_1/\sqrt{\pi})  + \frac{1}{2} \log(\log n^{\delta})} {\sqrt{2\delta\log n}}\,,
	& \mbox{if}\,~C_1 >0, \\
	\frac{\frac{1}{2}\log (C_2/(2\pi^2))}{\sqrt{2\delta\log n}}\,,
	& \mbox{if}\,~ C_1=0\,, \label{E2}
	\end{array}
	\right.                   
	\end{eqnarray} 
	with
	\begin{equation}\label{C1C2}
	C_1=\frac{1}{2c_K} \left(K^2(A)+K^2(-A) \right) ~~~~\mbox{and}~~~ C_2=\frac{1}{2c_K}\int [K'(t)]^2\, dt.~~~~~~~~~~~~~~
	\end{equation}
\end{thm}
The result in (\ref{D2}) can be used to construct confidence bands for $m(x)$. In fact, in light of (\ref{D2}), 
\begin{equation}\label{2**}
m_n(x)\pm \left(\frac{c_K\cdot \sigma^2_n(x)}{nh_n\cdot \widehat{f}_n(x)} \right)^{1/2}
\left(\frac{x^{(\alpha)}}{\sqrt{2\delta\log n}}+d_n\right), ~~ x\in[0,1]\,,
\end{equation}
represents an asymptotic $(1-\alpha)100$\% confidence band for the regression function $m(x)$  in the sense that,
as $n\to \infty$,
$$
P\left\{m(x)\in\,
m_n(x)\pm \left(\frac{c_K\cdot \sigma^2_n(x)}{nh_n\cdot \widehat{f}_n(x)} \right)^{1/2}
\left(\frac{x^{(\alpha)}}{\sqrt{2\delta\log n}}+d_n\right),~~ x\in[0,1]\right
\}
\longrightarrow 1-\alpha,
$$
where $x^{(\alpha)}=-(\log\log\left(\frac{1}{1- \alpha}\right) - \log 2)$ is the unique solution of the equation $\exp\{-2\exp(-x)\}=1-\alpha$. Furthermore, one can perform the test of hypothesis $H_0: m=m_0$, based on the statistic $t_n:=\sup_{x\in[0,1]}\,\big(\widehat{f}_n(x)/\sigma^2_n(x)\big)^{1/2}~\big| m_n(x)-m_0(x)\big|$, and reject $H_0$, at the significance level $\alpha$, if $t_n\, > ~ \sqrt{c_K/(nh_n)}\,\{(2\delta\log n)^{-1/2}
\left(\log 2-\log\log\left(\frac{1}{1-\alpha}\right)\right)+ d_n\}.$ \\

\subsection{The proposed estimator}
\vspace{2mm}\noindent
When the response variable $Y_i$ may be missing at random according to (\ref{MAR}), a very simple counterpart of (\ref{D1}) is the estimator that uses the complete cases only, i.e., the estimator
\begin{equation}\label{MBAR}
\overline{m}_n(x) \,:=\, \sum_{i=1}^n \Delta_i Y_i K((x-X_i)/h_n) \div  
\sum_{i=1}^n \Delta_i K((x-X_i)/h_n).
\end{equation}
The estimator in (\ref{MBAR}) is in a sense the right estimator. To appreciate this, observe that upon dividing the numerator and the denominator of the right hand side of (\ref{MBAR}) by the quantity $\sum_{i=1}^n K((x-X_i)/h_n)$, the estimator $\overline{m}_n(x)$ becomes the ratio of the kernel regression estimator of $E(\Delta Y|X=x)$ and the kernel regression estimator of $E(\Delta | X=x)$. Since 
$$\frac{E[(\Delta Y | X)}{E(\Delta | X)} \stackrel{\mbox{\tiny a.s.}}{=} \frac{E[E(\Delta Y | X,Y)|X]}{E(\Delta | X)} \stackrel{\mbox{\tiny by (\ref{MAR})}}{=} \frac{E[Y E(\Delta | X)|X]}{E(\Delta | X)} = E(Y|X),
$$
the estimator (\ref{MBAR}) is indeed a correct kernel-type regression estimator of $m(x)=E(Y|X=x)$. Despite its simplicity, there are no results available in the literature for the maximal deviations of $\overline{m}_n(x)$ similar to that in (\ref{D2}). Of course, it may be possible to establish results such as (\ref{D2}) for the maximal deviation of $\overline{m}_n(x)$  under the restrictive assumptions that $\sigma_0^2(x)$ or $p(x)=E[Y|X=x]$ are known, but in this paper we do not impose such assumptions. \\

\noindent
In what follows, we propose a kernel regression estimator of $m(x)$ where the presence of missing values is handled via a Horvitz-Thompson-type 
inverse weighting approach (Horvitz and Thompson (1952)). To motivate our proposed estimator, we first consider the simple but unrealistic case where the selection probability $p(x)=E[\Delta|X=x]$ is known. Now define
\begin{equation}\label{MTID}
\breve{m}_n(x)= \sum_{i=1}^n \, \frac{\Delta_i Y_i}{p(X_i)}\, K\left(\frac{x-X_i}{h_n}\right)\, \div  \,
\sum_{i=1}^n K\left(\frac{x-X_i}{h_n}\right),
\end{equation}
which is the kernel regression estimator of $E[\frac{\Delta Y}{p(X)}\big|X=x]= E[Y|X=x]$, where we have used the fact that 
\begin{equation}\label{MARnew}
E[\frac{\Delta Y}{p(X)}\big|X] \stackrel{\mbox{\tiny a.s.}}{=} \frac{1}{p(X)} E[E(\Delta Y | X,Y)|X] \stackrel{\mbox{\tiny by (\ref{MAR})}}{=} \frac{1}{p(X)} E[Y p(X) |X] \stackrel{\mbox{\tiny a.s.}}{=} E(Y|X).
\end{equation}
Therefore, when $p(x)$ is known, the estimator in (\ref{MTID}) is the kernel regression estimator of $E[Y|X=x]$. 
When $p(x)$ is unknown it can be replaced by an estimator $\widehat{p}_n(x)$; here we have in mind kernel regression estimators of $p(x)$. However, regardless of whether $p(x)$ or $\widehat{p}_n(x)$ is used in (\ref{MTID}), the maximal deviation of $\breve{m}_n(x)$ from $m(x)$ cannot be expected to yield the conclusion of Theorem \ref{THM-1}, as given by (\ref{D2}). This is because $\breve{m}_n(x)$ in (\ref{MTID})
is the kernel regression estimator of $E[Y^*|X=x]$, where $Y^*:= \frac{\Delta Y}{p(X)}$ does not always have a density with respect to the Lebesgue measure ($P\{Y^*=0\}\neq 0$ because $P\{\Delta=0\}>0$), which violates the assumption that the response variable ($Y^*$ in this case) has a pdf. To rectify this difficulty we start by artificially adding to $\frac{\Delta Y}{p(X)}$ a zero-mean continuous random variable $\epsilon$, whose pdf has a finite support, and where $\epsilon$ is independent of $(X, Y, \Delta)$. Clearly, the independence of $\epsilon$ and $X$ combined with the MAR assumption in (\ref{MAR}) yield $E\big[\frac{\Delta Y}{p(X)}+\epsilon\big|X=x\big]= E[Y|X=x]=m(x)$. The choice of the distribution of $\epsilon$ will be discussed later in Remark \ref{REM-1}. Now, let $\epsilon_1, \dots, \epsilon_n$ be iid copies of $\epsilon$, independent of the data $(X_i, Y_i, \Delta_i), ~i=1,\dots, n$, and consider the following revised version of (\ref{MTID}) 
\begin{equation}\label{MTIL}
\widetilde{m}_n(x)= \sum_{i=1}^n \, \left\{\left[\frac{\Delta_i Y_i}{p(X_i)}+\epsilon_i\right]\, K\left(\frac{x-X_i}{h_n}\right)\right\}\, \div  \,
\sum_{i=1}^n K\left(\frac{x-X_i}{h_n}\right),
\end{equation}
which is the kernel regression estimator of $m(x)$. Motivated by the naive estimator in (\ref{MTIL})
and the fact that $p(x)=E[\Delta|X=x]=E[\Delta+\epsilon | X=x]$, our proposed kernel-type regression estimator of $m(x)$ is given by
\begin{equation}\label{MHAT}
\widehat{m}_n(x)= \sum_{i=1}^n  \left\{\left[\frac{\Delta_i Y_i}{\widehat{p}_n(X_i)}+\epsilon_i\right]\, K\left(\frac{x-X_i}{h_n}\right)\right\}\, \div  \,
\sum_{i=1}^n K\left(\frac{x-X_i}{h_n}\right),
\end{equation}
where, for technical reasons that will be discussed under Remark \ref{REM-1}, we propose to use the following kernel-type estimator of $p(x)$ 
\begin{equation}\label{PHAT}
\widehat{p}_n(x)= \sum_{i=1}^n \left\{(\Delta_i +\epsilon_i)\, K\left(\frac{x-X_i}{\lambda_n}\right)\right\} \div  \,
\sum_{i=1}^n K\left(\frac{x-X_i}{\lambda_n}\right),
\end{equation}
instead of the usual kernel estimator $\widehat{p}_n(x)=\sum_{i=1}^n \Delta_i K((x-X_i)/\lambda_n) \div \sum_{i=1}^n K((x-X_i)/\lambda_n)$. Here, $\lambda_n\, (\neq h_n)$ is the smoothing parameter of the kernel.
Remark \ref{REM-1} below discusses the use of the random variable $\epsilon$,  its justification in the literature, and the choice of its distribution from both theoretical and applied points of view. Next, to establish the limiting distribution of the maximal deviation of $\widehat{m}_n(x)$ from $m(x)$, let 
\begin{equation}\label{WX}
Y^*=\frac{\Delta Y}{p(X)} + \epsilon~~~~\mbox{and put}~~~~~~
Z= Y^*- E[Y^*|X]~\Big(\stackrel{\mbox{\tiny a.s.}}{=} Y^*-m(X)\Big).   
\end{equation}
Let $\breve{G}(x,z)$ and $\breve{g}(x,z)$ be the joint cdf and the joint pdf of $(X,Z)$, respectively. Also, let $Q(z|x)$ and $q(z|x)$ be the conditional cdf and the conditional pdf of $Z$ given $X$ and consider the following counterpart of assumption (D):

\vspace{4mm}\noindent
{\bf Assumption (D$'$).} ~$\breve{g}^{1/2}(x,z)$ is differentiable with respect to both $x$ and $z$, and the partial derivatives are bounded. Furthermore, the inverse functions $Q^{-1}$ and $F^{-1}$ of $Q$ and $F$ exist and $\frac{\partial}{\partial x} Q^{-1}\big(z|F^{-1}(x)\big)$ and $\frac{\partial}{\partial z} Q^{-1}\big(z|F^{-1}(x)\big)$ are bounded.

\vspace{4mm}\noindent
{\bf Assumption (E$'$).} The kernel $K$ satisfies Assumption (E). Additionally, $K(x)=K(-x)$ and $K(x)\to 0$, as $|x|\to \infty$.

\vspace{4mm}\noindent
Regarding the selection probability $p(x)=E[\Delta|X=x]$ we assume

\vspace{2mm}\noindent
{\bf Assumption (F).} The selection probability $p(x)=P[\Delta=1|X=x]$ given by (\ref{MAR}) is twice differentiable with bounded derivatives. Also, $\inf_{x\in [a,b]}\, p(x) =:\,p_{0} > 0$, for some $p_0\in(0,1]$, where $[a,b]$ is as in Assumption (B). 

\vspace{4mm}\noindent
{\bf Assumption (G).} The random variables $\epsilon, \epsilon_1, \dots, \epsilon_n$ are iid zero-mean bounded random variables with a density function that vanishes off the interval $(a_0,b_0)$, for some $-\infty<a_0<b_0<\infty$. Also, $\epsilon_i$'s are independent of the data $(X_i, Y_i,\Delta_i)$.  

\vspace{4mm}\noindent
Here, the conditions $K(x)=K(-x)$ and $K(x)\to 0$, as $|x|\to \infty$ that appear under Assumption (E$'$) are as in Mack and Silverman (1982), whereas the second part of assumption (F) is standard in missing data literature and essentially amounts to requiring $Y$ to be observable with a non-zero probability for each $X=x$. Assumption (D$'$) is the counterpart of assumption (D). Next, let $\widehat{f}_n(x)$ be the kernel density estimator defined in (\ref{sig0}), and let 
\begin{equation}\label{SIGHAT}
\widehat{\sigma}^2_{\widehat{p}_n}(x)= \left\{\sum_{i=1}^n \, \left[\frac{\Delta_i Y_i}{\widehat{p}_n(X_i)}+\epsilon_i\right]^2\, K\left(\frac{x-X_i}{h_n}\right)\, \div  \,
\sum_{i=1}^n K\left(\frac{x-X_i}{h_n}\right)\right\} - \left(\widehat{m}_n(x)\right)^2\,,
\end{equation}
be the kernel regression estimator of the conditional variance $\sigma^2(x)=E\big[\big(\frac{\Delta Y}{p(X)}+\epsilon\big)^2\,\big|\,X=x\big]-\big(E\big[\frac{\Delta Y}{p(X)}+\epsilon\,\big|\,X=x\big]\big)^2~ \stackrel{\mbox{\tiny ($*$)}}{=} E\big[\big(\frac{\Delta Y}{p(X)}+\epsilon\big)^2\,\big|\,X=x\big]-m^2(x)$, where $(*)$ follows from (\ref{MAR}), (\ref{MARnew}), and assumption (G). Then we have the following result.
\begin{thm}\label{THM-2}
	Let $\widehat{m}_n(x)$, $\widehat{\sigma}^2_{\widehat{p}_n}(x)$, and $\widehat{f}_n(x)$ be as in (\ref{MHAT}), (\ref{SIGHAT}), (\ref{sig0}), respectively. Let $h_n=n^{-\delta}$ and $\lambda_n=n^{-\beta}$ for any $\delta$ and $\beta$ satisfying $1/5<\beta < \delta< 1/3$. Then, under assumptions (A)-(C), (D$'$), (E$'$), (F), and (G)
	\begin{equation*}
	P\left\{\sqrt{2\delta\log n}\left(
	\sqrt{\frac{nh_n}{c_K}} \,\, \sup_{x\in[0,1]}\,\sqrt{\frac{\widehat{f}_n(x)}{\widehat{\sigma}^2_{\widehat{p}_n}(x)}}~ \,\bigg| \widehat{m}_n(x)-m(x)\bigg|-d_n\right) \leq y\right\}\rightarrow  \exp\left(-2\,e^{-y}\right),~
	\end{equation*}
	as $n\to \infty$, where  $c_K=\int K^2(t)\,dt$ and $d_n$ is as in (\ref{E2}).
\end{thm}
This result immediately yields the following asymptotic $(1$$-$$\alpha)100$\% confidence bands for $m(x)$ when the response variable may be missing at random:
\begin{equation}\label{Cband}
\widehat{m}_n(x)\pm \left(\frac{c_K\cdot \widehat{\sigma}^2_{\widehat{p}_n}(x)}{nh_n\cdot \widehat{f}_n(x)} \right)^{1/2}
\left(\frac{x^{(\alpha)}}{\sqrt{2\delta\log n}}+d_n\right), ~~ x\in[0,1]\,,
\end{equation}
where $x^{(\alpha)}=-(\log\log\left(\frac{1}{1- \alpha}\right) - \log 2)$.

\vspace{5mm}
\begin{rem}\label{REM-1}
The main reason for employing the artificial variables $\epsilon_1, \dots,\epsilon_n$ in our methodology above is purely technical (and not quite necessary in numerical studies). The use of artificial or contrived variables in statistical estimation and inference is not new and, in fact, has a long history in the literature. Some classical examples include the problem of nearest neighbor classification when the $d$-dimensional covariate vectors do not have a pdf in which case the dimension is artificially increased to $d+1$ by including an additional random variable $\epsilon$ that has a pdf. This helps to establish the strong consistency of the nearest neighbor classifier and also works as a tie-breaking procedure (see, for example, Devroye et al. (1996, pp 175-176) or Gy\"orfi et al. (2002, p. 245)).  Another, perhaps more important, example of the use of artificial random variables is related to the weighted bootstrap approximation; see, for example, 
Mason and Newton (1992), Praestgaard and Wellner (1993), Janssen and Pauls (2003), Janssen (2005), Horv\'{a}th et al. (2000), Horv\'{a}th (2000), Burke (1998, 2000), Kojadinovic and Yan (2012), Kojadinovic, Yan, and Holmes (2011), and Mojirsheibani and Pouliot (2017) among others.   Since the conditional variance $\sigma^2(x)=E[(\frac{\Delta Y}{p(X)}+\epsilon)^2 | X=x] - m^2(x) = [p(X)]^{-1}E[Y^2|X=x] + E(\epsilon^2)-m^2(x)$, if $\epsilon$ is chosen to have a large variance, the estimator $\widehat{\sigma}^2_{\widehat{p}_n}(x)$ in (\ref{SIGHAT}) can be expected to be inflated. It therefore makes sense to choose $\epsilon$ to have a small variance. In fact, in our numerical work (Section \ref{Numerical}), we chose $\epsilon \sim$ Unif$\,[-a, a]$, where $a=10^{-3}$. However, our numerical results also show that one can actually replace $\epsilon_i$ by zero in (\ref{MHAT}) and (\ref{PHAT}), which is intuitively more appealing, and still expect to see the same numerical results. In other words, the presence of $\epsilon_i$ in (\ref{MHAT}) and (\ref{SIGHAT})  is only for theoretical purposes.
\end{rem}

\section{Numerical results}\label{Numerical}  In this section we carry out some simulation studies to evaluate (numerically) the finite-sample performance of the methods discussed in this paper. The results show that, in general, the proposed estimator performs well. We also take a close look at the performance of the complete-case estimator that is constructed based on the complete cases only. More specifically, in what follows we consider random samples $(X_1,Y_1),\dots,(X_n,Y_n)$ of sizes $n=200, 500,$ and $1000$ from the model
\[
Y= \sin\left(\pi[X^4+e^{\cos(X)}]\right) + \sigma(X) \cdot Z\,,~~ \mbox{with}~Z\sim N(0,1)\,,
\]
where $X\sim N(0.5 , 1)$ is independent of $Z$, and $\sigma^2(x)=1+e^{-(x+2)}$ represents the variance function $E(Y^2|X=x)-[E(Y|X=x)]^2$. 
Here, $Y_i$ can be missing at random based on one of the following two logistic missing probability  models for the function $\pi$ defined in (\ref{MAR}):

\vspace{2mm}
{\it Model A.} ~
$p(x) := P\{\Delta=1|X=x\}= \exp(1-2x)/\{1+\exp(1-2x)\}$.

\begin{figure}[th]
	\begin{center}
		\includegraphics[width=17.0cm, height=13.0cm]{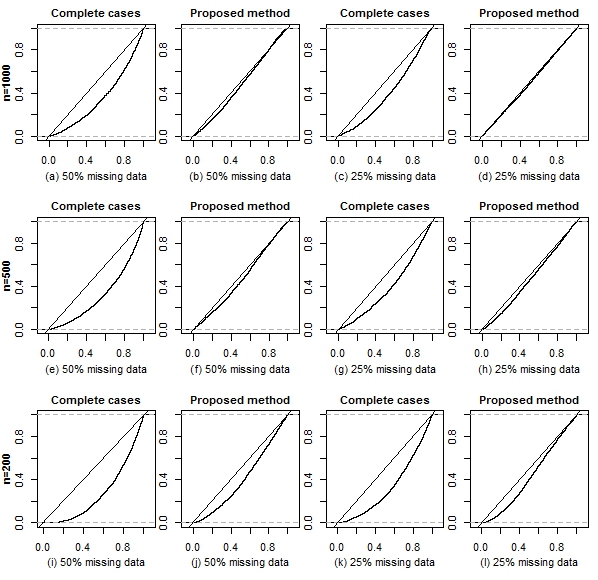}  
		\caption{\fontsize{10}{4}\selectfont Plots of empirical cdf's of $U_1, \dots, U_{3000}$ and $V_1,\dots, V_{3000}$ when $\epsilon_i \sim \mbox{Unif}(-\kappa, \kappa),~ \kappa=10^{-3}$ in (\ref{MHAT}) and (\ref{PHAT}).}  \label{FIG1}  
	\end{center}
\end{figure}

\vspace{3mm}
{\it Model B.} ~ $p(x) :=P\{\Delta=1|X=x\}=\exp(1+0.2x)/\{1+\exp(1+0.2x)\}$,

\vspace{2.5mm}\noindent
where $\Delta$=0 if $Y$ is missing (and $\Delta$=1, otherwise). These missing probability mechanisms yield roughly 50\% missing data under Model A and about 25\%  for Model B. Next, for our kernel estimators  $\widehat{p}_n(x)$ and $\widehat{m}_n(x)$ in  (\ref{PHAT}) and (\ref{MHAT}) and their smoothing parameters $h_n=n^{-\delta}$ and $\lambda_n=n^{-\beta}$
(subject to $1/5 < \beta < \delta < 1/3$), we employed the cross-validation approach of Racine and Li (2004) with the Epanechnikov kernel $K(u)=(3/4)(1-u^2)\cdot I\{|u|\leq 1\}$; this is implemented in the R package ``np'' (Hayfield and Racine 2008). 
As for the choice of  $\epsilon_i$'s that appear in (\ref{PHAT}) and (\ref{MHAT}), we considered $\epsilon_i \sim \mbox{Unif}(-\kappa\,,\,\kappa)$, $\kappa=10^{-3}$,  but we have also  considered the more appealing and practical choice of $\epsilon_i=0$, $i=1,\dots,n$.  Next, to evaluate the  performance of various estimators numerically, we computed the statistic
\begin{eqnarray}  \label{YnSTAT}
	U_n &:=& \sqrt{2\delta\log n}\bigg(
	\sqrt{nh_n/c_K} \, \sup_{x\in[0,1]}\,\sqrt{\widehat{f}_n(x)/\widehat{\sigma}^2_{\widehat{p}_n}(x)}~ \,\Big| \widehat{m}_n(x)-m(x)\Big|-d_n\bigg)
\end{eqnarray}
for sample sizes $n=200,\, 500, 1000$, and each of the missing probability models A and B, 
where $\widehat{m}_n(x)$, $\widehat{\sigma}^2_{\widehat{p}_n}(x)$, $\widehat{f}_n(x)$,  $c_K$, and $d_n$ are as in Theorem \ref{THM-2}.  In practice, to compute the supremum functional in (\ref{YnSTAT}), we used the maximum of $\sqrt{\widehat{f}_n(x)/\widehat{\sigma}^2_{\widehat{p}_n}(x)}~ \big| \widehat{m}_n(x)-m(x)\big|$ over a grid of 200 equally
spaced values of $x$ in the interval [0, 1]. Our initial pilot study shows that increasing the grid size to as large as 500 does not make any noticeable changes.  Next we note that if $n$ is not very small then by Theorem \ref{THM-2} the quantity 
$$
U= \exp\{-2\exp(-U_n)\}
$$ 
should be approximately a Unif\,[0,1] random variable. Repeating the entire above process a total of 3000 times yields $U_1, \dots,  U_{3000}$.  We also constructed the above statistic based on the complete cases only, i.e.,
\begin{eqnarray*} 
V_n &:=& \sqrt{2\delta\log n}\bigg(
\sqrt{nh_n/c_K} \, \sup_{x\in[0,1]}\,\sqrt{\bar{f}_n(x)/\bar{\sigma}^2_n(x)}~ \,\Big| \overline{m}_n(x)-m(x)\Big|-d_n\bigg),
\end{eqnarray*}
where $\overline{m}_n(x)$, $\bar{f}_n(x)$, $\bar{\sigma}^2_n(x)$ are the estimates of $m(x)$, $f(x)$, and $\sigma^2(x)$ based on the complete cases only; see (\ref{MBAR}) for the definition of the estimator $\overline{m}_n(x)$. If we put $V=\exp\{-2\exp(-V_n)\}$ then the 3000 Monte Carlo runs yield $V_1,\dots, V_{3000}$. Figure \ref{FIG1} gives plots of the empirical distribution functions of $U_1, \dots, U_{3000}$ and $V_1,\dots, V_{3000}$ for different sample sizes, different missing proportions, and the two choices of $\epsilon_i$'s. We have also included the $45^\circ$ line, which is the CDF of the Unif$\,[0,1]$ random variable. 

\begin{figure}[t]
	\begin{center}
		\includegraphics[width=17.0cm, height=13.0cm]{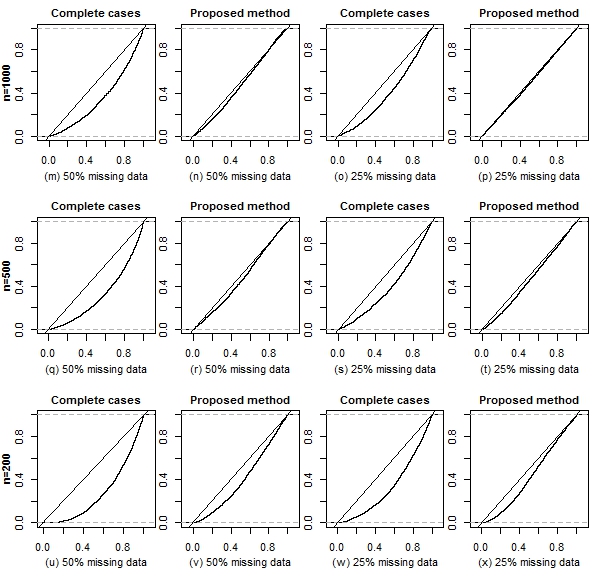}  
		\caption{\fontsize{10}{4}\selectfont Plots of empirical cdf's of $U_1, \dots, U_{3000}$ and $V_1,\dots, V_{3000}$ when $\epsilon_i=0$ in (\ref{MHAT}) and (\ref{PHAT}).}  \label{FIG2}  
	\end{center}
\end{figure}

Comparing plots (a) and (b) in Figure \ref{FIG1}, we see that the proposed estimator performs much better than the one based on complete cases when we have 50\% missing data and  $n=1000$; this is shown by the fact that the empirical CDF of $U_1, \dots, U_{3000}$ (which corresponds to the proposed estimator) is much closer to the $45^\circ$ line. In fact, as Figure \ref{FIG1} shows, the proposed estimator performs better at both 25\% and 50\% missing rates and for all sample sizes. Of course, the performance of the proposed estimator improves as $n$ increases, confirming the conclusion of Theorem \ref{THM-2}. From a practical point of view, the presence of the artificial random variables $\epsilon_1,\dots,\epsilon_n$ in the definition of the estimators (\ref{MHAT}) and (\ref{PHAT}) can be viewed as a nuisance. This is only needed as a technical tool and, in practice, one can take $\epsilon_i=0$. To confirm this, we also carried out the same simulation study with the more realistic choice of $\epsilon_i=0$, $i=1,\dots,n$ and the results were indistinguishable. Plots (m) to (x) in Figure \ref{FIG2} correspond to this setup. As Figure \ref{FIG2} shows, the results are virtually identical.
Next, we used our Monte Carlo simulation results to construct 90\% and 95\% uniform confidence bands for the regression function $m(x)$, over the set [0, 1]. This resulted in 3000 confidence bands for the regression function $m(x)$ for each sample size $n$ (= 200, 500, 1000) and each missing proportion (25\% and 50\%). The results summarized in Table \ref{Table-1} are for 90\% confidence bands. 
{
	\renewcommand{\arraystretch}{0.8} 
	\begin{table}[H]
		\caption{\fontsize{10}{5}\selectfont
			Coverage and the average area of confidence bands (averaged over 3000 bands) for the regression function $m(x)$; here $\kappa=10^{-3}$. Nominal coverage probability = 90\%.} 		
		\begin{center}
			\begin{tabular}{lr@{}c@{}c@{}c@{}c@{}c@{}c}
				\toprule
				& & \multicolumn{2}{c}{$n=200$}& \multicolumn{2}{c}{$n=500$} &\multicolumn{2}{c}{$n=1000$}\\ 
				\cmidrule(r){3-4}
				\cmidrule(l){5-6}
				\cmidrule(l){7-8}
				& Missing =& 50\% & 25\%   & 50\% & 25\% & 50\% & 25\%      \\ \hline\hline  
				Method      &   & & & & & &     \\ \hline  
				Via (\ref{Cband}), with   &Coverage =& 0.872 & 0.878 & 0.891 & 0.901 & 0.904 & 0.905   \\    
				$\epsilon_i\sim\,$Unif$(-\kappa, \kappa)$     &(Area) = &~1.884453&~\,1.440004~\,&1.416213~\,&1.150745 ~\,& 1.270975~\, & 0.820889\\ \hline      
				Via (\ref{Cband}), with  &Coverage =& 0.872 &0.878 &0.891 & 0.901 & 0.904 & 0.905 \\ 
				$\epsilon_1,\dots,\epsilon_n=0$  &(Area) = 
				&~1.884451&~\,1.440004~\,&1.416212~\,&1.150744 ~\,& 1.270973~\, & 0.820827\\ \hline
				Using complete    &Coverage =& 0.726 &0.794 &0.731 &0.836 & 0.778 & 0.859\\ 
				cases only  &(Area) =
				&~1.216772&~\,1.257136~\,&1.040574~\,&0.954291 ~\,& 0.841123~\, & 0.713148\\
				\bottomrule 
			\end{tabular}\label{Table-1}
		\end{center}
	\end{table}
\noindent
Here, coverage is computed as the proportion of confidence bands (out of 3000) that actually captured the true regression function $m(x)$ in the interval $[0, 1]$.  Table \ref{Table-1} also gives the average area of the 3000 confidence bands constructed under each setup. We make several observations here: (i) The table shows that, numerically, it makes no difference as to whether we consider
$\epsilon_i\sim\,$Unif$(-10^{-3}, 10^{-3})$ or the more intuitive and realistic  choice of $\epsilon_1,\dots,\epsilon_n =0$ in the estimators (\ref{MHAT}) and (\ref{PHAT}). This is clearly evident by the equal coverages in rows one and two of the Table \ref{Table-1} as well as the reported average areas (which are virtually the same up to 6 decimal places). As mentioned earlier, the presence of $\epsilon_i$'s in (\ref{MHAT}) and (\ref{PHAT}) are only for technical reasons. (ii) Table \ref{Table-1} also shows that the coverage of the bands based on complete cases can reduce far more noticeably as the missing rate changes from 50\% to 25\% than that of the proposed method. In other words, the proposed method is not as heavily influenced by the amount of missing data; this can be particularly important when a much larger proportion of the data is missing. (iii) We also note that the average areas of the bands based on our proposed method is somewhat higher than those based on complete cases. However, this does not mean that our bands are unnecessarily ``wider'' than what the theory suggests. 
In fact, as Theorem \ref{THM-2} as well as
Figures 1 and 2  show, the proposed bands are precisely those that are supported by the theory. 
Table \ref{Table-2} gives the corresponding results for  95\%  confidence bands. The conclusions are virtually the same those of Table \ref{Table-1}.

	\begin{table}[H]
		\caption{\fontsize{10}{5}\selectfont
			Coverage and the average area of confidence bands (averaged over 3000 bands) for the regression function $m(x)$; here $\kappa=10^{-3}$. Nominal coverage probability = 95\%.} 
		
		\begin{center}
			\begin{tabular}{lr@{}c@{}c@{}c@{}c@{}c@{}c}
				\toprule
				& & \multicolumn{2}{c}{$n=200$}& \multicolumn{2}{c}{$n=500$} &\multicolumn{2}{c}{$n=1000$}\\ 
				\cmidrule(r){3-4}
				\cmidrule(l){5-6}
				\cmidrule(l){7-8}
				& Missing =& 50\% & 25\%   & 50\% & 25\% & 50\% & 25\%      \\ \hline\hline  
				Method      &   & & & & & &     \\ \hline  
				Via (\ref{Cband}), with   &Coverage =& 0.932 & 0.938 & 0.947 & 0.957 & 0.949 & 0.954   \\    
				$\epsilon_i\sim\,$Unif$(-\kappa, \kappa)$     &(Area) = &~2.101318&~\,1.658559~\,&1.613458~\,&1.300233 ~\,& 1.423841~\, & 0.926563\\ \hline      
				Via (\ref{Cband}), with  &Coverage =& 0.932 &0.938 &0.947 & 0.957 & 0.949 & 0.954 \\ 
				$\epsilon_1,\dots,\epsilon_n=0$  &(Area) = 
				&~2.101317&~\,1.658558~\,&1.613456~\,&1.300233 ~\,& 1.423840~\, & 0.926561\\ \hline 
				Using complete    &Coverage =& 0.841 &0.887 &0.848 &0.918 & 0.872 & 0.928\\ 
				cases only  &(Area) =
				&~1.451007&~\,1.443447~\,&1.194201~\,&1.085454 ~\,& 0.955346~\, & 0.808382\\
				\bottomrule 
			\end{tabular}\label{Table-2}
		\end{center}
	\end{table}
}

\vspace{1mm}
\section{Concluding remarks} 
In this article, we have proposed a kernel-type method to construct asymptotically correct uniform confidence bands for an unknown regression function $m(x)=E[Y|X=x]$, over compact sets,  where the response variable $Y$ may be missing at random. The proposed method is fully nonparametric in that the selection probability,  the density function $f(x)$ of $X$, and the conditional variance of $Y$ given $X$ are all completely unknown. 
The proposed method is quite straightforward to implement and has good asymptotic properties. Furthermore, our numerical work shows that the proposed method has good finite-sample performance.  As explained in Remark \ref{REM-1}, the presence of the artificial variables, i.e.,  $\epsilon_i$'s, in our methodology and theoretical results are purely for technical reasons and, in practice, such variables can be taken to be zero.

\section{ Appendix: Proof of Theorem \ref{THM-2}}
We prove Theorem \ref{THM-2} in a number of steps.

\vspace{4mm}\noindent
\noindent\textbf{STEP 1.} Let $\widehat{m}_n(x)$ and $\widetilde{m}_n(x)$ be as in (\ref{MHAT}) and (\ref{MTIL}), respectively. Then, as $n\to\infty$, we have 
	\begin{equation}\label{L1}
	\sqrt{n h_n \log n}\, \sup_{x\in [0,1]} \,\left|\widehat{m}_n(x) - \widetilde{m}_n(x)\right| \longrightarrow ^p 0.
	\end{equation}
To show this, put $r_n(x)=\sum_{i=1}^n |\Delta_i Y_i|\, K((x-X_i)/h_n) / \sum_{i=1}^n K((x-X_i)/h_n)$ and observe that in view of Assumption (E)
\begin{eqnarray*}
\left|\widehat{m}_n(x) - \widetilde{m}_n(x)\right| &=&
\left|\sum_{i=1}^n \left[\Delta_i Y_i\,\left(\frac{1}{\widehat{p}_n( X_i)} -\frac{1}{p(X_i)}\right)  K\left(\frac{x-X_i}{h_n}\right)\right] \div \sum_{i=1}^n K\left(\frac{x-X_i}{h_n}\right)\right|\\
&\leq& \max_{1\leq i \leq n}\left[\left|\frac{\widehat{p}_n(X_i)-p( X_i)}{\widehat{p}_n(X_i)p(X_i)}\right|\, \mathbb{I}_{\big\{X_i\in [x-Ah_n,\,\,  x+Ah_n]\big\}}\right]\cdot r_n(x),
\end{eqnarray*}
where $\mathbb{I}_B$ denotes the indicator of a set $B$. 
But,  $r_n(x)\leq \max (|B_1|, |B_2|)$, by the second part of Assumption (A). Furthermore, in view of assumptions (A), (B), (C), (E$'$), and Theorem B of Mack and Silverman (1982) one has  $\sup_{x\in [0,1]} |\widehat{p}_n(x)-p(x)|=O_p\big(\sqrt{\log n/(n\lambda_n)}\,\big)$. Therefore, 
\begin{eqnarray}\nonumber\label{E21}
	\sup_{x\in [0,1]}\left|\widehat{m}_n(x) - \widetilde{m}_n(x)\right| &\leq&  \sup_{x\in [0,1]} \max_{1\leq i \leq n}\left[\left|\frac{\widehat{p}_n(X_i)-p( X_i)}{\widehat{p}_n(X_i)p(X_i)}\right| \mathbb{I}_{\big\{X_i\in[x-Ah_n,\,\, x+Ah_n]\big\}}\right]\cdot r_n(x)\nonumber\\
	&\leq& \frac{1}{p_0} \max_{1\leq i \leq n}\left|\frac{\widehat{p}_n(X_i)-p( X_i)}{\widehat{p}_n(X_i)} \right| \mathbb{I}_{\big\{X_i\in[-Ah_n,\,\, 1+Ah_n]\big\}}\cdot\sup_{x\in [0,1]}r_n(x)\nonumber\\
	&\leq& \frac{1}{p_0}\, \sup_{-A h_n\leq x\leq 1+Ah_n} \left|\frac{\widehat{p}_n(x)- p(x)}{\widehat{p}_n(x)}\right|\cdot\max (|B_1|, |B_2|)\nonumber \\
	&\leq& \frac{1}{p_0^2}\, O_p\left(\sqrt{\frac{\log n}{n\lambda_n}}\right)\cdot\max (|B_1|, |B_2|)  \, =\,O_p\left(\sqrt{\frac{\log n}{n\lambda_n}}\right)\\
	&& \mbox{(since $\lim_{n\to\infty}\mathbb{I}_{B_n}=\mathbb{I}_{\lim_{n\to\infty B_n}}$ for monotone sets $B_n$),~~~~~~~~\nonumber} 
\end{eqnarray}
where we have used the fact that $p_0 \leq \lim_{n\to\infty}\inf_{x}\, \widehat{p}_n(x) \leq 1$ (which  follows by noticing that $-\sup_{x} |\widehat{p}_n(x)-p(x)|+p_0 \leq \inf_x \widehat{p}_n(x) \leq \sup_x |\widehat{p}_n(x)-p(x)|+1$ and then taking the limit, as $n\to \infty$); here, the infimums are taken over the set $[-Ah_n \,, \, \mbox{1$+$$A$} h_n]$. Now  (\ref{L1}) follows from (\ref{E21}) together with the fact that $(\log n)^2 h_n/\lambda_n = (\log n)^2 \,n^{\beta-\delta}\to 0$, as $n\to \infty$ (because $\beta<\delta$).\\

\noindent \textbf{STEP 2.} Define the quantity
\begin{eqnarray}
\widetilde{\sigma}_n^2 &=& \left\{\sum_{i=1}^n \, \left[\frac{\Delta_i Y_i}{p(X_i)}+\epsilon_i\right]^2\, K\left(\frac{x-X_i}{h_n}\right)\, 
\div  \,\sum_{i=1}^n K\left(\frac{x-X_i}{h_n}\right)\right\} - \left(\widetilde{m}_n(x)\right)^2\,, ~~~~~~~~~\label{SIGTIL}
\end{eqnarray}
where $\widetilde{m}_n(x)$ is as in (\ref{MTIL}), and put  
\begin{eqnarray}
\sigma^2(x)&:=&
E\left[\left(\frac{\Delta Y}{p(X)}+\epsilon\right)^2 \bigg|\, X=x\right]-\left(E\left[\frac{\Delta Y}{p(X)}+ \epsilon \,\bigg|\, X=x \right]\right)^2 \nonumber \\
&=&
E\left[\left(\frac{\Delta Y}{p(X)}+\epsilon\right)^2 \bigg|\, X=x\right]-m^2(x), ~~(\mbox{by (\ref{MAR}), (\ref{MARnew}), and assumption (G))}.~~~~~~  \label{SIG}
\end{eqnarray} 
Also, let $\widehat{\sigma}^2_{\widehat{p}_n}(x)$  be as in (\ref{SIGHAT}). Then we have
\begin{eqnarray}\sup_{x\in [0,1]}\left|\widehat{\sigma}^2_{\widehat{p}_n}(x)-\widetilde{\sigma}^2_n(x)\right| &=& O_p\left(\sqrt{\frac{\log n}{n\lambda_n}}\,\right) + o_p\left(\frac{1}{\sqrt{n h_n \log n}}\right), \label{I}\\
\sup_{x\in [0,1]} \left|\widetilde{\sigma}^2_n(x)-\sigma^2(x)\right| &=& O_p\left(\sqrt{\frac{\log n}{n h_n}}\,\right),\label{II}
\end{eqnarray}
where $\sigma^2(x)$ is as in (\ref{SIG}). To establish (\ref{I}) and (\ref{II}), first observe that 
\begin{eqnarray*}
	\left|\widehat{\sigma}^2_{\widehat{p}_n}(x)-\widetilde{\sigma}^2_n(x)
	\right| &\leq& 
	\left| \left[\sum_{i=1}^n\left(
		\frac{1}{\widehat{p}^{\,2}_n(X_i)}-\frac{1}{p^2(X_i)}\right)\Delta_i Y_i^2 \cdot K\left(\frac{x-X_i}{h_n}\right)\right] \div\sum_{i=1}^n K\left(\frac{x-X_i}{h_n}\right)\right|\\
	&& + 2
	\left|\left[\sum_{i=1}^n\left(
		\frac{1}{\widehat{p}_n(X_i)}-\frac{1}{p(X_i)}\right)\epsilon_i \, \Delta_i Y_i\cdot K\left(\frac{x-X_i}{h_n}\right)\right]\div \sum_{i=1}^n K\left(\frac{x-X_i}{h_n}\right)\right|\\
	&& + \left|\Big(\widehat{m}_n(x)-\widetilde{m}_n(x)\Big)
	\Big(\widehat{m}_n(x)+\widetilde{m}_n(x)\Big)\right|	\\
	&:=& \big|V_{n,1}(x)\big|+ 2\big|V_{n,2}(x)\big| + \big|V_{n,3}(x)\big|
\end{eqnarray*}
But, by the second part of assumption (A), the second part of assumption (F), and the boundedness of the support of the distribution of $\epsilon$, we immediately find $\sup_{x\in [0,1]} |\widetilde{m}_n(x)|= O_p(1)$. Thus, by (\ref{L1}), $\sup_{x\in [0,1]} |\widehat{m}_n(x)| \leq \sup_{x\in [0,1]} |\widehat{m}_n(x)-\widetilde{m}_n(x)| + \sup_{x\in [0,1]} |\widetilde{m}_n(x)|= o_p\left((n h_n \log n)^{-1/2}\right)+O_p(1)=O_p(1)$. Therefore, by (\ref{L1}),
\[
\sup_{x\in [0,1]}\big|V_{n,3}(x)\big| \leq \sup_{x\in [0,1]} |\widehat{m}_n(x)-\widetilde{m}_n(x)|\cdot O_p(1) = o_p\left((n h_n \log n)^{-1/2}\right).
\]
Next, observe that 
\begin{eqnarray*}
\sup_{x\in[0,1]}\max_{1\leq i \leq n}\left|\frac{\widehat{p}^2_n(X_i)-p^2( X_i)} {\widehat{p}^2_n(X_i)p^2(X_i)}\right|\mathbb{I}_{\big\{X_i\in
[x-Ah_n,\,\, x+Ah_n]\big\}} &\leq& \frac{2}{p^2_0}\sup_{-A h_n\leq x\leq 1+A h_n}
\left|\frac{\widehat{p}_n(x)- p(x)}{\widehat{p}^2_n(x)}\right|\\
&=& O_p\big(\sqrt{\log n/(n\lambda_n)}\,\big),~~~~~~~~~~~~~~~
\end{eqnarray*}
which follows from the last part of assumption (F), the fact that $p_0^2 \leq \lim_{n\to\infty} \inf_{x} \,\widehat{p}^2(x) \leq \sup_x \,p^2(x) =1$, together with Theorem B of Mack and Silverman (1982), where, the infimum is taken over the set $[-Ah_n \,, \, \mbox{1$+$$A$} h_n]$. Therefore,
\begin{eqnarray*}
	\sup_{x\in [0,1]}\big|V_{n,1}(x)\big| &\leq& \sup_{x\in [0,1]}
	\left|\frac{\sum_{i=1}^n \Delta_i Y_i^2 \cdot K\left(\frac{x-X_i}{h_n}\right)}{\sum_{i=1}^n K\left(\frac{x-X_i}{h_n}\right)}\right| \times
	\frac{2}{p^2_0}\sup_{-A h_n\leq x\leq 1+A h_n}
	\left|\frac{\widehat{p}_n(x)- p(x)}{\widehat{p}^2_n(x)}\right|\\ 
	&=& O_p\big(\sqrt{\log n/(n\lambda_n)}\,\big),
\end{eqnarray*}
which follows because, in view of assumption (A), the first supremum term on the right side of the above inequality is bounded by $\max (B_1^2 , B_2^2)$. Similarly, we have $\sup_{x\in [0,1]}|V_{n,2}(x)| = O_p(\sqrt{\log n /(n\lambda_n)}\,)$, from which (\ref{I}) follows.
The proof of (\ref{II}) is rather straightforward  and goes as follows
\begin{eqnarray*}
	\sup_{x\in [0,1]} \left|\widetilde{\sigma}^2_n(x)-\sigma^2(x)\right| &\leq&
	\sup_{x\in [0,1]} \left| \frac{\sum_{i=1}^n \, \left[\frac{\Delta_i Y_i}{p(X_i)}+\epsilon_i\right]^2\, K\left(\frac{x-X_i}{h_n}\right)}{\sum_{i=1}^n K\left(\frac{x-X_i}{h_n}\right)} - E\left[\left(\frac{\Delta Y}{p(X)}+\epsilon\right)^2 \bigg|\, X=x\right]\right|\\
	&&+ \Big|\big(\widetilde{m}_n(x)-m(x)\big)
	\big(\widetilde{m}_n(x)+m(x)\big)\Big|\\
	&=&
	 O_p\left(\sqrt{\frac{\log n}{n h_n}}\,\right) + O_p\left(\sqrt{\frac{\log n}{n h_n}}\,\right) \,=\, O_p\left(\sqrt{\frac{\log n}{n h_n}}\,\right).
\end{eqnarray*}

\noindent \textbf{STEP 3.} Let $\widehat{f}_n(x)$, $\widehat{m}_n(x)$, $\widetilde{m}_n(x)$, $\widehat{\sigma}^2_{\widehat{p}_n}(x)$, and $\widetilde{\sigma}^2_n(x)$ be as in (\ref{sig0}), (\ref{MHAT}), (\ref{MTIL}), (\ref{SIGHAT}), and (\ref{SIGTIL}), respectively, and write
\begin{equation} \label{E01}
\sup_{x\in [0,1]} \sqrt{\frac{\widehat{f}_n(x)}{\widehat{\sigma}^2_{\widehat{p}_n}(x)}}\, \Big|\widehat{m}_n(x)-m(x)\Big| =
\sup_{x\in [0,1]} \sqrt{\frac{\widehat{f}_n(x)}{\widetilde{\sigma}^2_n(x)}}\, \Big|\widetilde{m}_n(x)-m(x)\Big| + R_n\,,
\end{equation}
where 
\begin{eqnarray}
	R_n &=& \sup_{x\in [0,1]} \sqrt{\frac{\widehat{f}_n(x)}{\widehat{\sigma}^2_{\widehat{p}_n}(x)}}\, \Big|\widehat{m}_n(x)-m(x)\Big| -
	\sup_{x\in [0,1]} \sqrt{\frac{\widehat{f}_n(x)}{\widetilde{\sigma}^2_n(x)}}\, \Big|\widetilde{m}_n(x)-m(x)\Big|   \nonumber \\
	&\leq& \sup_{x\in [0,1]} \sqrt{\frac{\widehat{f}_n(x)}{\widehat{\sigma}^2_{\widehat{p}_n}(x)}}\, \Big|\widehat{m}_n(x)-\widetilde{m}(x)\Big|  \nonumber \\
	&& ~~~~~~~~ +
	\left(\sup_{x\in [0,1]} \sqrt{\frac{\widetilde{\sigma}^2_n(x)}{\widehat{\sigma}^2_{\widehat{p}_n}(x)}}-1\right) \cdot 
	\sup_{x\in [0,1]} \sqrt{\frac{\widehat{f}_n(x)}{\widetilde{\sigma}^2_n(x)}}\, \Big|\widetilde{m}_n(x)-m(x)\Big|    \nonumber \\
	&:=& R_n(1) + R_n(2).   \label{E1}
\end{eqnarray}
To deal with the supremum on the r.h.s of (\ref{E01}), we note that  $\widetilde{m}_n(x)$ and $\widetilde{\sigma}_n^2(x)$ that appear in this supremum term are, respectively, the kernel regression estimator of $E(Y^*|X=x)$ and the kernel estimator of the conditional variance of $Y^*$, as given by (\ref{SIG}), based on the iid ``data'' $(X_i, Y^*_i),$ $i=1,\dots,n$, where $Y^*$ is given by (\ref{WX}). 
It is straightforward to see that when assumptions (A) and (F) hold then $P\{B^*_1\leq Y^* \leq B^*_2\}=1$, where $B_1^*=\min(0,B_1)+a_0$~ and $B_2^*=\frac{B_2}{p_0}+b_0$, with the constants $B_1$ and $B_2$ as in assumption (A), and where $a_0$ and $b_0$ are as in assumption (G). Also, in view of assumptions  (A) and (G), the random vector $(X, Y^*)$ has a pdf. Therefore,  when assumption (A) holds for the distribution of $(X,Y)$ then, 
because of asumption (F),  it also holds for the distribution of $(X,Y^*)$. Similarly,  if $\sigma_0^2(x)$ satisfies assumption (C) then so does $\sigma^2(x)$ (in view of assumption (F)); to show this, simply observe that in view of (\ref{MAR}) we have $\sigma^2(x)\stackrel{\mbox{\tiny via (\ref{SIG})}}{=}[(p(x))^{-1}-1] E(Y^2|X=x)+E(\epsilon^2)+\sigma_0^2(x)$.
Therefore, as a consequence of Theorem \ref{THM-1}, under assumptions (A), (B), (C), (D$'$), (E), (F), and (G),
\begin{equation}\label{EE2}
\sqrt{2 \delta \log n}\,\left\{ \sqrt{\frac{n h_n}{c_K}}\, \sup_{x\in [0,1]} \sqrt{\frac{\widehat{f}_n(x)}{\widetilde{\sigma}^2_n(x)}}\, \Big|\widetilde{m}_n(x)-m(x)\Big|- d_n\right\} \longrightarrow^{d}\, Y\,,
\end{equation}
where $P(Y\leq y)=\exp\left\{-2 \exp(-y)\right\}$, $y\in \mathbb{R}$, $c_K=\int K^2(t)\,dt$, and $d_n$ is as in (\ref{E2}). 
Therefore to prove Theorem \ref{THM-2}, it is sufficient to show that $\sqrt{n h_n \log n}\,\big(R_n(1)+R_n(2)\big)\to^p 0$, as $n\to\infty$. First we show that $\sqrt{n h_n \log n}\,|R_n(2)|\to^p 0$. To show this, observe that by (\ref{EE2}) 
\begin{equation}\label{EE3}
\sup_{x\in [0,1]} \sqrt{\frac{\widehat{f}_n(x)}{\widetilde{\sigma}^2_n(x)}}\, \Big|\widetilde{m}_n(x)-m(x)\Big| = O_p\left(\sqrt{\log n/(n h_n)}\right).
\end{equation}
Furthermore 
$
\big|\sup_x\sqrt{\widetilde{\sigma}^2_n(x)/\widehat{\sigma}^2_{\widehat{p}_n}(x)}-1\big|
\leq \sup_{x}\big|\sqrt{\widetilde{\sigma}^2_n(x)/\widehat{\sigma}^2_{\widehat{p}_n}(x)}-
1\big| \leq \frac{\sup_{x}\big|\widetilde{\sigma}^2_n(x)-
	\widehat{\sigma}^2_{\widehat{p}_n}(x)\big|}{\inf_{x} \widehat{\sigma}^2_{\widehat{p}_n}(x)}.$
But by (\ref{I}), $\sup_{x\in [0,1]}\left|\widehat{\sigma}^2_{\widehat{p}_n}(x)-\widetilde{\sigma}^2_n(x)\right| = O_p\left(\sqrt{\log n/(n \lambda_n)}\right) + o_p\left((n h_n \log n)^{-1/2}\right)$. We also note that 
$\,\inf_x\, \widehat{\sigma}^2_{\widehat{p}_n}(x) \geq \inf_x\, \left\{\widehat{\sigma}^2_{\widehat{p}_n}(x)-\sigma^2(x)\right\} +\inf_x\, \sigma^2(x) \geq - \sup_x \big|\widehat{\sigma}^2_{\widehat{p}_n}(x)-\sigma^2(x)\big| +\inf_x\,\sigma^2(x)
$, where $\sigma^2(x)$ is as in (\ref{SIG}). Similarly, observe that $\inf_x\,\widehat{\sigma}_n^2(x) \leq \sup_x \big|\widehat{\sigma}^2_{\widehat{p}_n}(x)-\sigma^2(x)\big| +\sup_x\,\sigma^2(x)$. Thus we have
\[
- \sup_x \big|\widehat{\sigma}^2_{\widehat{p}_n}(x)-\sigma^2(x)\big| +\inf_x\,\sigma^2(x) \,\leq
\inf_x\,\widehat{\sigma}_n^2(x) \leq\, \sup_x \big|\widehat{\sigma}^2_{\widehat{p}_n}(x)-\sigma^2(x)\big| +\sup_x\,\sigma^2(x).
\]
Now, in view of (\ref{I}) and (\ref{II}), and upon taking the limit in the above chain of inequalities, as $n\to\infty$,  we find $0 < \lim_{n\to\infty} \inf_x \, \widehat{\sigma}_n^2(x) <\infty$, which yields 
$$
\bigg|\sup_{x\in [0,1]}\sqrt{\widetilde{\sigma}^2_n(x)/\widehat{\sigma}^2_{\widehat{p}_n}(x)}-1\bigg|
= O_p\left(\sqrt{\log n/(n \lambda_n)}\right) + o_p\left((n h_n \log n)^{-1/2}\right).
$$
This in conjunction with (\ref{EE3}) imply that 
$\sqrt{n h_n \log n}\,|R_n(2)|\to^p 0$, where $R_n(2)$ is as in (\ref{E1}). Next, observe that $\sup_x \big|\widehat{f}_n(x)/\widehat{\sigma}_n^2(x)\big| \leq \big[\sup_x \big|\widehat{f}_n(x)-f(x)\big|+\sup_x f(x)\big]/\inf_x \widehat{\sigma}^2_{\widehat{p}_n}(x) = O_p(1)$, which follows because $\sup_x|\widehat{f}_n(x)-f(x)|=o_p(1)$ and by the fact that
 $0 < \lim_{n\to\infty} \inf_x \, \widehat{\sigma}_n^2(x) <\infty$ (as shown above). Combining these results, we have
 \begin{equation*}
 \sqrt{n h_n \log n}\,|R_n(1)| \leq \,
\sqrt{\sup_{x\in [0,1]}\left|\widehat{f}_n(x)/\widehat{\sigma}^2_{\widehat{p}_n}(x)\right|}\,\cdot \sup_{x\in [0,1]} \Big|\widehat{m}_n(x)-\widetilde{m}(x)\Big|  
= O_p(1) \cdot o_p(1) \, = \, o_p(1). \label{E2B}
 \end{equation*}
 Putting the above results together, we have $\sqrt{n h_n \log n}\,|R_n|=o_p(1)$.
Theorem \ref{THM-2} now follows from this together with   (\ref{E01}), (\ref{E1}), and (\ref{EE2}).

\hfill $\Box$

\vspace{20mm}
\noindent {\Large \bf References}
\vspace{2.3mm}
 
\noindent \vspace{2.3mm}\noindent Burke, M.: A Gaussian bootstrap approach to estimation and tests in Asymptotic Methods in Probability and Statistics. E. (eds.) B. Szyszkowicz, pp. 697-706. North-Holland, Amsterdam (1998) \\

\vspace{2.3mm}\noindent Burke, M.: Multivariate tests-of-fit and uniform confidence bands using a weighted bootstrap. Statist. Probab. Lett. 46, 13-20 (2000)\\

\vspace{2.3mm}\noindent Cai, T., Low, M., Zongming, M.: Adaptive confidence bands for nonparametric regression functions. J. Amer. Statist. Assoc. 109, 1054-1070 (2014)  \\

\vspace{2.3mm}\noindent Claeskens, G., Van Keilegom, I.: Bootstrap confidence bands for regression curves and their derivatives. Ann. Statist. 31, 1852-1884 (2003)\\

\vspace{2.3mm}\noindent Devroye, L., Gy\"{o}rfi, L., Lugosi, G.: A probabilistic theory of pattern recognition. Springer-Verlag, New York (1996) \\

\vspace{2.3mm}\noindent Eubank, R.L., Speckman, P.L.: Confidence Bands in Nonparametric Regression. J. Amer. Statist. Assoc. 88 (424), 1287-1301 (2012) \\

\vspace{2.3mm}\noindent Gu, L., Yang, L.: Oracally efficient estimation for single-index link function with simultaneous confidence band. Electron. J. Stat. 9, 1540-1561 (2015) \\

\vspace{2.3mm}\noindent Gy\"orfi, L., Kohler, M., Krzy\.zak, A., Walk, H.: A distribution-free theory of nonparametric regression. Springer-Verlag, New York (2002)\\

\vspace{2.3mm}\noindent H\"ardle, W.: Asymptotic maximal deviation of M-smoothers. J. Multivariate Anal. 29, 163-179 (1989)\\ 

\vspace{2.3mm}\noindent H\"ardle, W.: Applied Nonparametric Regression. Cambridge University Press (1990)\\

\vspace{2.3mm}\noindent H\"ardle, W., Song, S.: Confidence bands in quantile regression. Econometric Theory 26 (4), 1-22 (2010) \\

\vspace{2.3mm}\noindent Hollander, M., McKeague., I.W., Yang, J.: Likelihood ratio-based confidence bands for survival functions. J. Amer. Statist. Assoc. 92, 215-227 (1997) \\

\vspace{2.3mm}\noindent Horv\'{a}th, L. Approximations for hybrids of empirical and partial sums processes. J. Statist. Plann. Inference. 88:1-18 (2000)\\

\vspace{2.3mm}\noindent Horv\'{a}th, L., Kokoszka, P., Steinebach, J.: Approximations for weighted bootstrap processes with an application. Statist. Probab. Lett. 48, 59-70 (2000) \\

\vspace{2.3mm}\noindent Horvitz, D.G., Thompson D.J. A generalization of sampling without replacement from a finite universe. J. Amer. Statist. Assoc. 47, 663-685 (1952) \\

\vspace{2.3mm}\noindent Janssen, A.: Resampling Student's t-type statistics. Ann. Inst. Statist. Math. 57, 507-529 (2005)\\

\vspace{2.3mm}\noindent Janssen, A., Pauls, T.: How do bootstrap and permutation tests work? Ann. Statist. 31, 768-806 (2003)\\

\vspace{2.3mm}\noindent Johnston, G.J.: Probabilities of maximal deviations for nonparametric regression function estimates. J. Multivariate Anal. 12, 402-414 (1982) \\

\vspace{2.3mm}\noindent Kojadinovic, I., Yan., J.: Goodness-of-fit testing based on a weighted bootstrap: A fast large-sample alternative to the parametric bootstrap. Canad. J. Statist. 40, 480-500 (2012) \\

\vspace{2.3mm}\noindent Kojadinovic, I., Yan, J., Holmes, M.: Fast large-sample goodness-of-fit for copulas. Statist. Sinica 21, 841-871 (2011) \\

\vspace{2.3mm}\noindent Konakov, V.D., Piterbarg, V.I. : On the convergence rate of maximal deviation distribution. J. Multivariate Anal. 15, 279-294 (1984)\\

\vspace{2.3mm}\noindent Liero, H.: On the maximal deviation of the kernel regression function estimate. Series Statistics 13, 171-182 (1982)  \\

\vspace{2.3mm}\noindent Lei, Q., Qin, Y.: Confidence intervals for nonparametric regression functions with missing data: multiple design case. J. Syst. Sci. Complex. 24, 1204-1217 (2011) \\

\vspace{2.3mm}\noindent Little, R.J.A., Rubin, D.B: Statistical analysis with missing data. Wiley, New York (2002) \\

\vspace{2.3mm}\noindent Li, G., van Keilegom, I.: Likelihood ratio confidence bands in nonparametric regression with censored data, Scand. J. Statist. 2, 547-562 (2002)\\

\vspace{2.3mm}\noindent Mack, Y.P., Silverman, Z.: Weak and strong uniform consistency of kernel regression estimates. Z. Wahrsch. Verw. Gebiete 61, 405-415 (1982) \\

\vspace{2.3mm}\noindent Mason, D.M., Newton, M.A.: A rank statistics approach to the consistency of a general bootstrap. Ann. Statist. 20, 1611-1624 (1992) \\

\vspace{2.3mm}\noindent Mass\'{e}, P., Meiniel, W.: Adaptive confidence bands in the nonparametric fixed design regression model. J. Nonparametr. Stat. 26, 451-469 (2014)  \\

\vspace{2.3mm}\noindent Mojirsheibani, M., Pouliot, W.: Weighted bootstrapped kernel density estimators in two sample problems.J. Nonparametr. Stat. 29, 61-84 (2017)\\

\vspace{2.3mm}\noindent Mondal, S., Subramanian, S.: Simultaneous confidence bands for Cox regression from semiparametric random censorship. Lifetime Data Anal. 22, 122-144 (2016) \\

\vspace{2.3mm}\noindent Muminov, M.S.: On the limit distribution of the maximum deviation of the empirical distribution density and the regression function. I. Theory Probab. Appl. 55, 509-517  (2011) \\

\vspace{2.3mm}\noindent Muminov, M.S.: On the limit distribution of the maximum deviation of the empirical distribution density and the regression function II. Theory Probab. Appl. 56, 155-166  (2012) \\

\vspace{2.3mm}\noindent Nadaraya, E.A.: Remarks on nonparametric estimates for density functions and regression curves. Theory Probab. Appl. 15, 134-137 (1970)\\

\vspace{2.3mm}\noindent Neumann, M.H., Polzehl, J.: Simultaneous bootstrap confidence bands in nonparametric regression. J. Nonparametr. Stat. 9, 307-333 (1998) \\

\vspace{2.3mm}\noindent Praestgaard, J., Wellner, J.A.: Exchangeably weighted bootstraps of the general empirical process. Ann. of Probab. 21, 2053-2086 (1993)\\

\vspace{2.3mm}\noindent Proksch, K.: On confidence bands for multivariate nonparametric regression. Ann. Inst. Statist. Math. 68, 209-236 (2016)\\

\vspace{2.3mm}\noindent Qin, Y., Qiu, T., Lei, Q.: Confidence intervals for nonparametric regression functions with missing data. Comm. Statist. Theory Methods. 43, 4123-4142 (2014)\\

\vspace{2.3mm}\noindent Racine, J., Li, Q.: Cross-validated local linear nonparametric regression. Statist. Sinica 14, 485-512 (2004) \\

\vspace{2.3mm}\noindent Hayfield, T., Racine, J.: Nonparametric econometrics: The np package. Journal of statistical software (2008) \\

\vspace{2.3mm}\noindent Song, S., Ritov, Y., H\"ardle, W.: Bootstrap confidence bands and partial linear quantile regression. J. Multivariate Anal. 107, 244-262 (2012)\\

\vspace{2.3mm}\noindent Wandl, H.: On kernel estimation of regression functions. Wissenschaftliche Sitzungen zur Stochastik, WSS-03, Berlin. (1980)\\

\vspace{2.3mm}\noindent Wang, Q., Shen, J.: Estimation and confidence bands of a conditional survival function with censoring indicators missing at random. J. Multivariate Anal. 99, 928-948  (2008)\\

\vspace{2.3mm}\noindent Wang, Q., Qin, Y.: Empirical likelihood confidence bands for distribution functions with missing responses. J. Statist. Plann. Inference 140, 2778-2789 (2010)\\

\vspace{2.3mm}\noindent Watson, G.S.:  Smooth regression analysis. Sankhya Ser. A 26, 359-372 (1964)\\

\vspace{2.3mm}\noindent Xia, Y.: Bias-corrected confidence bands in nonparametric regression.  J. R. Stat. Soc. Ser. B. Stat. Methodol. 60, 797-811 (1998)

\end{document}